\documentclass[12pt]{amsart}

 \newcommand{\new}{\newcommand}                        
 \new{\trunc}{{\otimes}}                               
 \new{\tnsr}{\otimes}                                  
 \new{\tensor}{\otimes}                                
 \new{\iso}{\cong}                                     
 \new{\union}{\cup}                                    
 \new{\W}{\mathfrak{W}}                                
 \new{\g}{\mathfrak{g}}                                
 \new{\h}{\mathfrak{h}}                                
 \new{\uqg}{U_q(\g)}                                   
 \new{\bracket}[1]{\langle#1\rangle}                   
 \new{\qdim}{\operatorname{qdim}}                      
 \new{\Hom}{\operatorname{Hom}}                        
 \new{\C}{\mathbb{C}}                                     
 \new{\R}{\mathbb{R}}                                     
 \new{\Cat}{\mathfrak{C}}                              
 \new{\ZZ}{\mathfrak{Z}}                               
 \new{\F}{\mathfrak{F}}                                
 \new{\Z}{\mathbb{Z}}                                     
 \new{\Q}{\mathbb{Q}}                                     
 \newcounter{letter}                                   
 \newenvironment{alist}{
 \begin{list}{{(\alph{letter})}}{\usecounter{letter}}
 }{\end{list}}                                         
 \newtheorem{thm}{Theorem}                             
 \newtheorem{prop}{Proposition}                        
 \newtheorem{lem}{Lemma}                               
\theoremstyle{definition}                              

\begin{document}
\title{Closed subsets of the Weyl alcove and TQFTs}
\author{Stephen F. Sawin}
\address{Math  Dept.\\  Fairfield U.\\(203)254-4000x2573\\
ssawin@fair1.fairfield.edu}
\begin{abstract}
For an arbitrary simple Lie algebra $\g$ and an arbitrary root of unity
$q,$ the  closed subsets of the Weyl alcove of the quantum group
$U_q(\g)$ are classified.   Here a closed subset is a set such that
if any two weights in the Weyl alcove are in the set,  so is any
weight in the Weyl alcove which corresponds to an irreducible summand of
the tensor product of a pair of representations with highest weights the two
original weights. The ribbon category associated to each closed subset admits a ``quotient'' by a trivial subcategory as described by Brugui\`eres (\cite{Bruguieres00}) and M\"uger (\cite{Muger99}), to give a modular category and a framed three-manifold invariant or a spin modular category and a spin three-manifold invariant (\cite{Sawin02b}).  Most of these theories are equivalent to theories defined in (\cite{Sawin02a,Sawin02b}), but several exceptional cases represent the first nontrivial examples to the author's knowledge of theories which contain noninvertible trivial objects, making the theory much richer and more complex.
\end{abstract}
\maketitle

\section*{Introduction}

Quantum groups, i.e.  quantized universal enveloping algebras
of simple Lie algebras, together with their representation theory have
been the subject of much fruitful investigation, and are of interest
from many perspectives, but one particularly important application is
to link and three-manifold invariants.  The general setting is that
these quantum groups are ribbon Hopf algebras, and hence their
representation theory forms a ribbon category, from which one can
construct an invariant of links and more generally labeled graphs
embedded in $S^3$ with similar properties to the original Jones
polynomial.   When the complex parameter $q$ on which these algebras depend is a
root of unity, their representation theory satisfies the more
restrictive requirements of a  \emph{modular category,} from which one can
construct a three-manifold invariant which satisfies Atiyah's axioms
for topological quantum field theory.  More specifically, the set
of  representations spanned by  the subset of the irreducible
representations in the Weyl alcove, with the ordinary tensor
product of representations replaced by the truncated tensor
product, forms a modular category. 

Any subset of the Weyl alcove which spans a set of representations
closed under both duality and the truncated tensor product determines a new
ribbon category which is  a full subcategory of the original.  Of course this
subcategory encodes only a subset of the link information of the original theory,
and on this basis does not seem of interest.  However, if this ribbon
subcategory happens to be modular, there is no reason to think that
the resulting three-manifold invariant is determined by the original
three-manifold invariant,
and in fact it is not apparent that it has any connection with the
original.   Thus finding closed subsets of the Weyl alcove which are
modular is an important question.   

In fact, the modularity requirement can be relaxed quite a bit.
M\"uger \cite{Muger99} and Brugui\`eres \cite{Bruguieres00} have shown
that under favorable circumstances which hold for the quantum group
examples (e.g., the existence of a unitary structure) a ribbon
category admits a kind of quotient which yields a modular category in
the absence of a certain  easily identifiable obstruction.   In fact, even in
the presence of this obstruction \cite{Sawin02b} a similar process
yields an invariant of spin three-manifolds.  Thus a classification of
the closed subsets of the Weyl alcove, together with an identification
of the quotient and when the quotient yields a modular category, would
give a complete summary of the invariants of three-manifolds which can
be constructed out of the Weyl alcove in this fashion.  The present
article classifies the closed subsets of the Weyl alcove and describes 
the subcategory of so called degenerate objects by which one quotients.

A second reason for considering the closed subsets of the Weyl alcove
is that they might plausibly correspond to quotients of the quantum
group, or at least of a subalgebra.  In fact in the classical case the
closed subsets of the Weyl 
chamber correspond exactly to quotients of the simply-connected groups
by a subgroup of the center: i.e., there is a one-to-one
correspondence between closed subsets of the Weyl chamber on the one
hand and Lie groups with the given 
Lie algebra on the other.   Thus quotients associated to closed subsets of the Weyl
alcove, if they exist, might be viewed as quantum analogues of the
nonsimply-connected groups.

Finally, in the course of the classification we shall construct some 
unexpected theories at level $k=2.$  Some of these theories 
(associated to the quantum group of type $B_n$ and $D_n$)  appear to give new 
invariants and admit skein relations that suggest we might be able to 
compute for these theories much of what can be computed in the 
$\rm{SU}(2)$  (i.e., $A_1$) theory.  Nevertheless, these theories 
exhibit novel behavior (specifically,  the 
subcategory of degenerate objects is the representation category of a 
nonabelian group) worthy of further study.    

The analysis of the closed subsets of the Weyl alcove was begun by the
author in
\cite{Sawin02a}.  There closed subsets which correspond precisely to
the classical closed subsets, i.e. to nonsimply-connected Lie groups,
were found and their invariants identified.   These subsets correspond to
Chern-Simons theories for nonsimply-connected Lie groups conjectured
by Dijkgraaf and Witten \cite{DW90}.  A second
collection of closed subsets was also identified and classified there.
These subsets, associated with certain corners of the Weyl alcove,   formed
very simple ribbon categories  (in fact the category associated to the
group algebra of an abelian group, with a mildly deformed $R$-matrix)
and manifold invariants depending only on homology which had been
studied by Murakami, Ohtsuki and Okada \cite{MOO92}.  The present article completes the
identification of the closed subsets of the Weyl alcove, demonstrating
that the closed subsets given above are the only ones except for certain special
cases at level $k=2.$

The paper is organized as follows.  Section 1  introduces the
needed facts about Lie algebras, quantum groups and the truncated
tensor product, which can all be found in Humphreys
\cite{Humphreys72}, Kirillov \cite{Kirillov96}, Kassel
\cite{Kassel95}, Turaev \cite{Turaev94} and \cite{Sawin02a}.  Sections
2 and 3 give the classification of the closed subsets, and Section 4
identifies when the quotient is modular and what the resulting
invariant is for the exceptional cases not discussed in
\cite{Sawin02a}.   Finally Section 5 shows that the closure under
duality condition on closed subsets is implied by the closure under
tensor products.  This observation is fairly independent of the rest
of the paper, but is a natural question.  In particular it is of
great relevance when using skein theory and cabling to explore link
and three-manifold invariants (for example work of Turaev and Wenzl
\cite{TW93} and of Wenzl \cite{Wenzl93}).

\subsection*{Acknowledgments}  I would like to thank I. Frenkel and
A. Liakhovskaia for helpful conversations and suggestions.

\section{Quantum Groups and the Weyl Alcove}

Let $\g$ be a simple Lie algebra   
and let 
$\{\alpha_i\}_{i \leq r}$ be the simple roots of $\g.$ The weight
lattice $\Lambda$ is spanned by  the \emph{fundamental weights}
$\{\lambda_i\}_{i\leq r}$  given by $(\lambda_i,\alpha_j)=
\delta_{i,j}(\alpha_i,\alpha_i)/2.$ 

  The Weyl group is
denoted by $\W,$ and the set of weights in the fundamental Weyl chamber is
called
$\Lambda^+$  (we will loosely refer to this set itself as the Weyl chamber).  
Half the sum of the positive roots is called $\rho$, the unique short
root in the Weyl chamber is called $\beta,$ and  the unique
long root in the Weyl chamber is called $\theta$ (in the simply-laced
case either will refer to the unique root in the Weyl alcove).
The root $\theta$   is the  highest weight of the adjoint
representation of $\g.$  The dual Coxeter 
number $\check{h}$ is defined to be $(\rho,\theta)+1,$ the value of the quadratic
Casimir on the adjoint representation.

Let $q=e^{2\pi i/(k+\check{h})},$ for some natural number $k$.   Recall that
there is an irreducible representation  $V_\lambda$  of the quantum
group $U_q(\g)$ for each weight $\lambda$ in the \emph{Weyl alcove}
$\Lambda_0,$ $i.e.,$ weights $\lambda \in \Lambda^+$ such that
$(\lambda,\theta) \leq k.$  
Kirillov \cite{Kirillov96} shows that the category of representations of the quantum
group which are a direct sum of representations in the Weyl alcove forms
a  semisimple ribbon $*$-category if the ordinary tensor product is
replaced by the 
\emph{truncated tensor product,} $\trunc,$ 
which is  the maximal subspace of the ordinary tensor product 
 isomorphic to a direct sum of representations
 in the Weyl alcove. Kirillov's $q$ is normalized a bit differently from what we use
here, which follows the conventions of \cite{Sawin02a}, but the previous
sentence is the essential thing we take from Kirillov, and determines
the normalization.  The truncated tensor product operation  on the 
lattice of isomorphism classes of representations in the category 
(with direct sum as addition) forms
a commutative, distributive, associative 
multiplication with the trivial
representation $V_0$ acting as identity,    determined by
\begin{equation} V_\lambda \trunc V_\gamma \iso
\bigoplus_{\eta \in \Lambda_0} N_{\lambda,
\gamma}^\eta \odot
V_\eta,\end{equation}
where $N_{\lambda, \gamma}^\eta$ are nonnegative integers  and $N
\odot V$ indicates the direct sum of $N$ copies of the representation $V.$
  The author (\cite{Sawin05b})  gives the following  formula for these numbers (generalizing a result of Andersen and Paradowski (\cite{AP95}))
\begin{equation}  \label{eq:AP}
N_{\lambda,\gamma}^\eta = \sum_{\sigma \in
\W_0} (-1)^\sigma m_\gamma(\lambda-\sigma(\eta))
\end{equation}
where $m_\lambda(\mu)$ is the dimension of the $\mu$ weight space inside the
classical representation of highest weight $\lambda$ and $\W_0$ is the
quantum Weyl group, which is 
generated by  reflection about the hyperplanes
$\{x|(x+\rho,\alpha_i)=0\}$ for each simple root  $\alpha_i$ together with $\{x| 
(x,\theta)=k+1\}.$ Also
\begin{equation}\label{eq:theta}
C_\lambda=q^{(\lambda,\lambda+2\rho)/2},\end{equation} 
and $\qdim(\lambda)>0,$ where $\qdim(\lambda)$ is the invariant of the
unknot and $C_\lambda$ is the factor a full twist applies to the link invariant.

For simplicity, since we will only ever need to consider
representations up to isomorphism, we will confuse the (isomorphism
class of the) representation $V_\lambda$ with the weight $\lambda,$
for example writing $\lambda \trunc \gamma= \bigoplus_\eta
N_{\lambda,\gamma}^\eta \odot \eta.$  Caution should be used with the
operations $\oplus$ and $\odot,$ because the weights are elements of
the weight lattice and therefore admit a lattice  addition and scalar
multiplication denoted by $\lambda+\gamma$ and $n\lambda,$
respectively, which we will also make frequent use of.  Note that
$\lambda+\gamma \neq \lambda \oplus \gamma,$ and $n\lambda \neq n
\odot \lambda.$   In particular, $0$ is the additive identity in the
lattice, but is the multiplicative identity for $\trunc.$   We hope
that the brevity of the notation outweighs this modest awkwardness.

The following  theorem is proven in \cite{Sawin02a}.
\addtocounter{thm}{-1}
\begin{thm}   \label{th:extreme}
There is an injection $\ell$ from $Z(G),$ the center of the simply connected Lie group with Lie algebra $\g,$  to the fundamental weights of the 
Weyl alcove such that $z$ acts on the classical representation 
$V_\gamma$ as $\exp(2 \pi i (\gamma,\ell(z))) \cdot \operatorname{id}_\gamma.$  The 
fundamental weights $\lambda_i$ in the image of $\ell$ are exactly those 
for which $(\lambda_i,\theta)=1$ and the associated root $\alpha_i$ is 
long, and are also exactly those for which there is a unique element $\tau_i$ 
of the classical Weyl group sending the standard base to the base 
$\{\alpha_j\}_{j \neq i} \cup \{-\theta\}.$ If we 
define $\phi_i(\gamma)=k\lambda_i +\tau_i(\gamma),$ then $\phi_i$ is an 
isometry of the Weyl alcove and of the simplex 
$\{\lambda:(\lambda,\alpha_j) \geq 0 \text{ and } (\lambda,\theta)\leq k\}$ 
and $\phi_i(\lambda \trunc \gamma)=\phi_i(\lambda) \trunc \gamma$   If
we use $k$ also to  
represent the map on the weight lattice which multiplies each weight 
by the 
number $k,$ then $k\ell$ is a homomorphism in the sense that 
$k\ell(zz')=k\ell(z) \tensor k\ell(z').$  Weights in the range of $k\ell$ can be 
characterized as extreme points of the simplex  $\{\lambda:(\lambda,\alpha_j) 
\geq 0 \text{ and } (\lambda,\theta)\leq k\}$ such that a neighborhood of 
the weight $0$ intersected with the simplex is isometric to a neighborhood 
of the extreme point intersected with the alcove, the isometry being given 
by $\phi_i.$
\end{thm} 

If $Z$ is a subgroup of $Z(G),$ let $\Delta_{Z}$ be the image of $Z$ 
under $k\ell.$    The
subset  of the Weyl chamber consisting of weights $\gamma$ such that
$Z$ acts trivially on  $V_\gamma$ is the
intersection of the Weyl chamber with a sublattice of the weight lattice, and
its elements are in one-to-one correspondence with representations of the Lie
group $G/Z.$  The intersection $\Gamma_Z$ of this set with the Weyl alcove may
be thought of loosely as the ``Weyl alcove for quantized $G/Z,$'' and consists
of those weights $\gamma$ in the alcove for which $(\gamma,\ell(z))$ is an
integer for all $z \in \Z.$   

A few key facts about the truncated tensor product were proven in
\cite{Sawin02a}.

\begin{lem} \label{lm:exist} 
For any $\sigma$ in the classical Weyl group $\W,$ and any weights
$\gamma,$ $\lambda$ in the Weyl alcove, if
$\lambda+\sigma(\gamma)$ is in the Weyl alcove,
 then
$\lambda \trunc \gamma$ contains $\lambda+\sigma(\gamma)$ as a summand.
\end{lem}

\begin{lem} \label{lm:theta}
$\lambda \trunc \lambda^\dagger$ contains as a summand $\theta$ if $k \geq
2$ and $\lambda$ is not a corner  (i.e. a multiple of a fundamental
weight such that $(\lambda,\theta)=k$).  In the nonsimply-laced case it
contains as a summand $\beta$ unless $(\lambda,\alpha_i)=0$ for every
short simple root $\alpha_i.$ 
\end{lem}

We say a fundamental weight $\lambda_i$ is \emph{long} or \emph{short}
according to whether $\alpha_i$ is long or short.  We say a
\emph{long} weight $\lambda_i$
is  \emph{sharp} or \emph{dull} according to whether
$\bracket{\lambda_i,\theta}=1$ or not.  Finally, we say that a weight
$\lambda$ is a \emph{long, short, sharp or dull corner} if
$\bracket{\lambda,\theta}=k$ and it is a multiple of a fundamental
weight with the corresponding property.  Thus for example the above
lemmas tell us that $\lambda \tensor \lambda^\dagger$ contains
$\theta$ or $\beta$ as a summand unless $\lambda$ is a dull corner.

\section{Closed Subsets}

Following \cite{Sawin02a}, we say that a subset $\Gamma$ of the set of
representations in the Weyl alcove $\Lambda_0$ is \emph{closed} if it is closed
under duals (for every $\gamma \in \Gamma,$ the weight $\gamma^\dagger$ of the dual
representation  is in $\Gamma$) and under the truncated tensor
product (if $\lambda, \gamma \in \Gamma$ then every $\eta$ such that
$N_{\lambda,\gamma}^\eta \neq 0,$ i.e. such that  $\eta$ corresponds
to  a direct summand of  $\lambda \trunc \gamma,$ is in $\Gamma$). We will see in 
Section \ref{se:tensor} that in the quantum groups case the second condition
implies the first.   Such subsets correspond exactly to ribbon subcategories,
and following M\"uger \cite{Muger99} if they meet a certain easily checked
condition they admit a quotient which is modular and gives a TQFT and
three-manifold invariant.   The following theorem gives a complete
classification of closed subsets of the Weyl alcove.  The proof encompasses
this and the next section.

\begin{thm} \label{th:closed}
The closed subsets of the Weyl alcove are the following
\begin{alist}
\item For any subgroup $Z \subset Z(G),$ the set $\Gamma_Z.$
\item For any subgroup $Z \subset Z(G),$ the set $\Delta_Z.$ 
\item  For $k=2,$ 
\begin{enumerate}
\item   $E_7$ the set
$\{0,\lambda_6\},$ 
\item $E_7$ with set $\{0,\lambda_2,2\lambda_7\},$
\item   $E_8$  with set $\{0,\lambda_1\},$
\item   $B_l$ with set $\{0,2\lambda_1,\lambda_j,\lambda_{2j}, \ldots,
\lambda_{(n-1)j/2}\}$ where $j>2$ and $2l+1=nj,$
\item $D_{l}$  with set $\{0, 2\lambda_1, \lambda_j,\lambda_{2j},
\ldots, \lambda_{(n-1)j},2\lambda_{l-1},2\lambda_l\}$ where $j>2$ and
$l=nj,$ 
\item $D_l$ with set $\{0,
2\lambda_1, \lambda_j,\lambda_{2j}, 
\ldots, \lambda_{(n-1)j/2}\}$ where $j>2$ and $2l=nj$ for $n$ odd.
\end{enumerate}
  Here $\lambda_i$ are the fundamental weights ordered as in
\cite{Humphreys72}. We will call (1)-(6)   the exceptional closed
subsets.
\end{alist}
\end{thm}

Before embarking on a proof of the theorem, some general discussion
is in order.    Formula (\ref{eq:AP}) is difficult to use in general,
but we will find it suffices for most of our purposes to examine it carefully
only when one of the factors is $\theta$ or $\beta.$ Our strategy will
be to show that for any $\lambda$ which is not in one of these
exceptional cases and not a sharp corner(and therefore  not in the image of $k  \ell$) an appropriate tensor
product of factors of $\lambda$ and $\lambda^\dagger$ contains $\theta$ or
$\beta.$ It will be easy to see that any closed subset
containing $\theta$ or $\beta$ is of the first type in the theorem.  By
Lemma \ref{lm:theta},   $\theta$ or $\beta$ is contained
in a tensor product of copies of $\lambda$ and $\lambda^\dagger$
unless  $\lambda$ is a long corner.  

By Lemma \ref{lm:exist}, $\lambda \trunc \beta$ will contain as a
summand with multiplicity one every weight $\lambda+\alpha$ for which $\alpha$
is a root (or a short root in the nonsimply laced case) and $\lambda +
\alpha$ is in the Weyl alcove.  Thus a crucial question is this : Given any
dull corner $\lambda,$  for which 
$\alpha$ is $\lambda+\alpha$ in the Weyl alcove?  Of course if
$\lambda_i+\alpha$ is in the Weyl alcove for $k=(\theta,\lambda_i)$
(i.e. when $\lambda_i$ is a corner),
then $n\lambda_i + \alpha$ is in the Weyl alcove whenever $n\lambda_i$
is a corner.

Below we list for each dull $\lambda_i$ the roots (or in the nonsimply-laced case short roots) $\alpha$ such
that $\lambda_i+\alpha$ is in the Weyl 
alcove for $k=(\lambda_i,\theta).$  The labeling of the roots and fundamental
weights are as in
\cite{Humphreys72}.   The roots are written out as a sum of simple
roots and $\theta$ in such a 
way that each additional simple root has negative inner product with
the sum up to that point (as can be checked using the Cartan matrix
(\cite{Humphreys72}, page 59), confirming recursively that the sum is a
root.  That the root is short in the nonsimply-laced case and that the
entire sum has nonnegative inner product with each simple root and
inner product with $\theta$ less than or equal to that of $\lambda_i$
can be checked from the Cartan matrix and the expansion of $\theta$
(\cite{Humphreys72}, page 66).

\begin{itemize}
\label{chart}
\item[$B_l$:]  For $2 \leq i \leq l-1.$ $\lambda_i + \alpha_{i+1} +
\alpha_{i+2} + 
\cdots + \alpha_l,$ $\lambda_i-(\alpha_i+\alpha_{i+1} + \cdots +\alpha_l).$

\item[$C_l$:]  $\lambda_l+\alpha_1 +\alpha_2 + \cdots +
\alpha_{l-1},$  $\lambda_l-(\alpha_1+ \cdots + \alpha_l).$

\item[$D_l$:]  For $2\leq i \leq l-2.$  $\lambda_i-(\alpha_{i-1} + \alpha_i
+ \cdots + \alpha_{l-2} + \alpha_{l-1} + \alpha_l +\alpha_{l-2} +
\cdots +\alpha_i),$
  $\lambda_i+\alpha_{i+1} + \alpha_{i+2} + \cdots + \alpha_{l-2} +
\alpha_{l-1}+\alpha_{l} + \alpha_{l-2} + \alpha_{l-3} + \cdots +
\alpha_{i+2}$  (for 
the cases $i=l-2,$ and  $i=l-3$ the last formula should read $\lambda_{l-2} +
\alpha_{l-1} $ and $\lambda_{l-3} +\alpha_{l-2} + \alpha_{l-1} + \alpha_l$
respectively).

\item[$E_6$:]  $\lambda_2+\alpha_1+\alpha_3 + \alpha_4+\alpha_5+\alpha_6,$
$\lambda_2-\theta.$  $\lambda_3+\alpha_1,$
$\lambda_3-\theta+\alpha_2+\alpha_4+\alpha_5+\alpha_6,$ $\lambda_4+\alpha_1+\alpha_3,$
$\lambda_4+\alpha_5+\alpha_6,$ $\lambda_4-\theta+\alpha_2.$
$\lambda_5+\alpha_6,$ $\lambda_5-\theta+\alpha_2 + \alpha_4 + \alpha_3+\alpha_1.$

\item[$E_7$:]
$\lambda_1 + \alpha_3 + \alpha_4 + \alpha_5
+\alpha_6 + \alpha_7 +\alpha_2 + \alpha_4 +\alpha_5 +\alpha_6,$
$\lambda_1-\theta.$  $\lambda_2 
-\theta+\alpha_1+\alpha_3 + \alpha_4 + \alpha_5 +\alpha_6 + \alpha_7.$
$\lambda_3+\alpha_2+\alpha_4+\alpha_5+\alpha_6 +\alpha_7,$
$\lambda_3-\theta -\alpha_1.$
$\lambda_4+\alpha_2,$ $\lambda_4+\alpha_5+\alpha_6+\alpha_7.$
$\lambda_5 +\alpha_6+\alpha_7,$
$\lambda_5-\theta +\alpha_1+\alpha_3+\alpha_4+\alpha_2.$
$\lambda_6 +\alpha_7,$
$\lambda_6-\theta + \alpha_1 + \alpha_3 + \alpha_4+\alpha_5+\alpha_4+\alpha_3.$

\item[$E_8$:]
$\lambda_1-\theta+\alpha_8 + \alpha_7 + \alpha_6 + \alpha_5 + \alpha_4
+ \alpha_3 + \alpha_2 + \alpha_4 + \alpha_5 +\alpha_6 + \alpha_7,$
$\lambda_1-(\alpha_3+ \alpha_4 + \alpha_2 + \alpha_5 + 
\alpha_6 + \alpha_7 + \alpha_8 + \alpha_ 4 + \alpha_5 + \alpha_6 +
\alpha_7).$ 
$\lambda_2 -\theta +\alpha_1+ \alpha_3 + \alpha_4 +\alpha_5 + \alpha_6
+\alpha_7 +\alpha_8,$ $\lambda_2-(\alpha_2+\alpha_4+\alpha_3+\alpha_1+
\alpha_5+\alpha_6 + \alpha_4  +\alpha_3 +\alpha_5 +
\alpha_4 +\alpha_2).$  
$\lambda_3+\alpha_1,$
$\lambda_3-(\alpha_1+\alpha_3+ \alpha_4 +\alpha_2 + \alpha_5 +\alpha_4
+ \alpha_3).$  
$\lambda_4+\alpha_1+\alpha_3,$ $\lambda_4 + \alpha_2.$
$\lambda_5+ \alpha_1+\alpha_3+\alpha_4+\alpha_2,$  
$\lambda_5
-\theta +  \alpha_6 + \alpha_7 +
\alpha_8.$  $\lambda_6 + \alpha_1+\alpha_3+ \alpha_4 + \alpha_2,$  
$\lambda_6 - \theta +  \alpha_7 +  \alpha_8.$  
$\lambda_7 + \alpha_1+\alpha_3 + \alpha_4 +
\alpha_2 + \alpha_5 + \alpha_6 + \alpha_4 + \alpha_3  +
\alpha_5 + \alpha_4 + \alpha_2,$ 
$\lambda_7 - \theta+ \alpha_8.$
$\lambda_8 + \theta-\alpha_8 - \alpha_7 -\alpha_6 - \alpha_5 -
\alpha_4 -\alpha_3 -\alpha_2 -\alpha_4 -\alpha_5 -\alpha_6 - \alpha_7 -\alpha_8,$
$\lambda_8-\theta.$

\item[$F_4$:]
$\lambda_1+ \alpha_2+\alpha_3 + \alpha_4 + \alpha_3,$
$\lambda_1-(\alpha_1+ \alpha_2+\alpha_3).$  $\lambda_2 + \alpha_3 +
\alpha_4,$  $\lambda_2 - (\alpha_2 + \alpha_3)$.

\item[$G_2$:]  $\lambda_2+\alpha_1,$  $\lambda_2- (\alpha_1 +
\alpha_2).$
\end{itemize}

The first observation to be made from this list is that for each of
the nonsimply-laced algebras except $B_l$ there is a short root $\alpha$ such that
$\theta+ \alpha$ is in the Weyl alcove for $k \geq 2$ (For  $C_l$
$\theta=\lambda_l$ and $\alpha=\alpha_1 + \alpha_2 + \cdots +
\alpha_{l-1},$ for $F_4$ $\theta=\lambda_1$ and
$\alpha=\alpha_2+\alpha_3+\alpha_4+\alpha_3,$ and for $G_2$
$\theta=\lambda_2$ and $\alpha=\alpha_1$).  Thus by Formula (\ref{eq:AP})
$N_{\theta,\theta}^{\theta+\alpha}$ contains a contribution from $\sigma=1$ since
$m_\theta(\alpha)=1.$  In order for it to contain a contribution for some other
$\sigma,$  that $\sigma$ would have to be a reflection about a short root $\alpha_i$
such that $(\theta+\alpha,\alpha_i)=0$ and $\alpha-\alpha_i$ is long.   If
$\alpha-\alpha_i$ is long then $(\alpha,\alpha_i)=0.$  One can easily
check by direct computation that there is no such $\alpha_i$
 in any of these case.  Therefore $\theta \trunc \theta$ contains
$\theta+\alpha,$ so $\theta \trunc \theta \trunc \theta$ contains $(\theta + \alpha)
\trunc \theta$ which by Lemma \ref{lm:exist} contains $\beta.$   For
$B_l$ when $k>2$ the same argument applies to $\theta+\beta.$  When
$k=2$ we will see below that a power of $\theta$ contains
$2\lambda_l,$  which is a short corner and hence a higher power
contains $\beta.$  Thus we conclude that if a closed subset of the Weyl alcove
contains $\theta,$ it also contains $\beta.$

\begin{lem} \label{lm:coset}
If a closed subset of the Weyl alcove contains  $\beta,$ it
is  of the form $\Gamma_Z$ for some $Z.$
\end{lem}

\begin{proof}
Consider $\lambda$ in the root lattice and in the Weyl alcove, and
choose a path
in the root lattice connecting $\lambda$ to $0$ such that the difference between
successive points in the path is a short root.    By reflecting
about hyperplanes of reflection in the quantum Weyl group we can replace it by such
a path crossing fewer such hyperplanes, and by induction can find such a path
entirely within the alcove.  Then by Lemma \ref{lm:exist} $\lambda$ is contained
in $\beta^{\trunc n},$ where $n$ is the length of this path.  Thus if a closed
subset contains  $\beta$  it
contains all of the root lattice $\Lambda_r.$

If $\lambda, \gamma$ are in the Weyl alcove and in the same coset of
$\Lambda/\Lambda_r,$ their difference is in the root lattice, and thus any closed
subset containing one of these and $\beta$ contains the other.  So any closed subset
containing $\beta$ is a union of cosets intersected with the Weyl alcove.  The
tensor product of two weights is a nonempty sum (because the quantum dimensions are
nonzero) and is in the product of their cosets, so the set of cosets making up such
a closed subset is a subgroup of $\Lambda/\Lambda_r.$  This proves that the closed subset
is in the intersection of the preimage of a subgroup of $\Lambda/\Lambda_r$ with the Weyl
alcove.   By Theorem \ref{th:extreme} the map sending  $z \in
Z(G)$  to  $\exp(2\pi i(\ell(z),\cdot)) \in (\Lambda/\Lambda_{r})^{*}$ is a
 group isomorphism.     So such a
subgroup is dual to some subgroup $Z\subset Z(G),$ and the closed subset is  exactly $\Gamma_Z.$
\end{proof}

Now if $\lambda$ is a dull corner it is a multiple of one of the weights in the chart
above. If $\lambda_i +\alpha$ is in the Weyl alcove for
$k=(\lambda_i,\theta),$ then  $n\lambda_i+\alpha$ is in the Weyl
alcove for $k=n(\lambda_i,\theta).$  So  from the chart there are at
least two elements of the Weyl alcove a short root away from $\lambda,$ except if
$\lambda$ is a multiple of  $\lambda_2$ for $E_7.$   In this
exceptional case
 $\lambda_2/2-\alpha_2$ is 
in the Weyl alcove, so there
are still at least two elements of the alcove a short root away from $\lambda$
 except when $\lambda=\lambda_2$ and $k=2.$
Thus with this exception  
if $\lambda$ is a dull corner  then by Lemma \ref{lm:exist} $\lambda
\trunc \beta$ contains two weights.   But each of these weights when tensored with
$\beta$ contain $\lambda,$ so $\lambda \trunc \beta \trunc \beta$ contains $\lambda$
with multiplicity at least $2.$ Therefore $\lambda \trunc \lambda^\dagger$ contains two
weights in $\beta \trunc \beta.$  Of course one is the trivial weight, so the other
must be of the form $\beta + \alpha$ for $\alpha$ a short root different from
$-\beta.$  If this $\beta + \alpha$ is not a corner dual to a long
root, a tensor power of it then contains $\beta.$

\section{Proof of the Theorem}

\begin{proof}[Proof of Theorem \ref{th:closed}]
It was shown in
\cite{Sawin02a} that the subsets of the form $\Gamma_Z$ and $\Delta_Z$ are closed .   We will see below that the exceptional cases are closed by computing
the truncated tensor product completely.  So we have only to show that every closed subset
is of this form.

We will assume the closed set contains some $\lambda$  not in the image of $k  \ell,$
and show it is either the exceptional case or it contains $\theta$ or $\beta,$ in
which case by Lemma \ref{lm:coset}  it falls into the first category.  By Lemma
\ref{lm:theta} we may assume $\lambda$ is a corner dual to a long root.

\subsection{If
$\mathbf{k=1}$} There are no dull corners, so there is nothing to prove.

\subsection{If $\mathbf{k=2}$} Except for $\lambda_2$ of $E_7$ if $\lambda$ is a
corner dual to a long root and not in the range of $k  \ell$ then
$\lambda \trunc \lambda^\dagger$ 
contains something nontrivial in 
$\beta \trunc \beta,$  so we may assume that $\lambda$ is a such a
corner and $\lambda \trunc \lambda^\dagger$ contains a summand of the
form $\beta + \alpha$ for $\alpha$ a short root.

\begin{itemize}

\item For $A_l,$ there are no dull corners.

\item For $B_l,$ there is nothing to prove if $l=2,$ so assume $l>2.$
The Weyl alcove consists of $\lambda_i$ for $i\leq l$ and $2\lambda_1$ and
$2\lambda_l.$  By checking which differences among these are short
roots and noting $\lambda_1=\beta$ we conclude 
\begin{align*}
\lambda_i \trunc \lambda_1 &=
\begin{cases} \lambda_{i-1}\oplus  \lambda_{i+1} & \text{for $1<i<l-1$}\\
 0 \oplus 2\lambda_1 \oplus \lambda_2   & \text{for $i=1$}\\
 2 \lambda_l \oplus \lambda_{l-2} & \text{for $i=l-1$}\\
\lambda_l & \text{for $i=l,$}
\end{cases}\\
2\lambda_l \trunc \lambda_1 &= 2\lambda_l \oplus \lambda_{l-1},\\
2\lambda_1 \trunc \lambda_1 &=\lambda_1,\\
2\lambda_1 \trunc 2\lambda_1 &=0.
\end{align*}
From this we conclude recursively
\begin{align*}
2\lambda_1 \trunc \lambda_i &= \lambda_i \qquad \text{for $i\leq l$},\\
2\lambda_1 \trunc 2\lambda_l&= 2\lambda_l,\\
\lambda_i \trunc \lambda_j&=
\begin{cases}\lambda_{i-j} \oplus \lambda_{i+j} & \text{for $l>i >
        j$ and $i+j <l$}\\
        \lambda_{i-j} \oplus 2\lambda_{l} & \text{for $l > i >
        j$ and $i+j = l, l+1$}\\
        \lambda_{i-j} \oplus \lambda_{2l+1-i-j} &\text{for $l>i>j$ and $i+j>l+1$}\\
         0 \oplus 2\lambda_1 \oplus \lambda_{2i}
        & \text{for $l>i=j$ and $2i<l$}\\
         0 \oplus 2\lambda_1 \oplus 2\lambda_{l}
        & \text{for $l>i=j$ and $2i = l, l+1$}\\
        0 \oplus 2\lambda_1 \oplus \lambda_{2l+1-2i} & \text{for $l>i=j$ and $2i>l+1.$}
\end{cases}
\end{align*}

Notice first of all that a closed subset containing $\lambda_2=\theta$
will contain $\lambda_i$ for all even $i<l.$ In particular it must
contain either $\lambda_{l-2}$ or $\lambda_{l-1},$ so it must contain
$2\lambda_l,$ and therefore since this is a short corner it
must contain $\beta.$  Thus we recover our promised assertion that
even for $B_l$ at level $k=2,$ a closed subset containing $\theta$
contains $\beta$ and therefore is of the from $\Gamma_Z.$

Let $\Gamma$ be a closed subset which is not of the form $\Gamma_Z$
or $\Delta_Z.$  We know that $\Gamma$ cannot contain
$\lambda_1,$ $\lambda_2,$ $\lambda_l$ or $2\lambda_l,$ or else it
would be of the form $\Gamma_Z,$ and $\Gamma$ must contain something
other than $0$ and $2\lambda_1.$  So let $j$ be the least $j$ such
that $\lambda_j \in \Gamma.$  Necessarily $l>j>2.$  By the product
rules above $\lambda_j, \lambda_{2j}, \lambda_{3j},\ldots \in \Gamma_Z.$
Suppose $m$ is the largest such that $\lambda_{mj} \in \Gamma.$  Then
every summand of $\lambda_{mj} \trunc \lambda_j$ is in $\Gamma.$
Clearly $(m+1)j \geq l,$ and in fact $(m+1)j>l+1,$ or else $2\lambda_l
\in \Gamma.$  So we conclude that $\lambda_{2l+1-(m+1)j} \in \Gamma.$
Now $mj<l<(m+1)j$ so $2l+1-(m+1)j$ is within $j$ of $mj.$  If they are
not equal then $\lambda_{mj} \trunc \lambda_{2l+1-(m+1)j}$ contains
$\lambda_i$ where $i$ is this difference, contradicting the minimality
of $j.$  Thus we conclude $mj=2l+1-(m+1)j,$ and $\Gamma$ contains
 set (4) in the statement of the theorem, with $n=2m+1.$  If it 
contained any $\lambda_{i}$ not in this set, there would be $p$ such 
that $|i-pj|<j,$ and hence $\lambda_{|i-pj|}$ would be in the set, 
contradicting the minimality of $j.$

\item For $C_l$ there is nothing to prove since  $\theta$ is the only
      dull corner.

\item For $D_l$ at $k=2$ the weights are $\lambda_i$ for $1 \leq i
\leq l,$ $2\lambda_1,$ $2\lambda_{l-1},$  $2\lambda_l,$ and
$\lambda_{l-1} + \lambda_l.$  By checking which differences among
these are roots we see
$$\lambda_i \trunc \theta = 
\begin{cases}
\lambda_1 \oplus \lambda_3 & \text{for $i=1$}\\
0 \oplus 2\lambda_1 \oplus \lambda_4 & \text{for $i=2$}\\
\lambda_{i-2} \oplus \lambda_{i+2} &
\text{for $2<i<l-3$}\\
\lambda_{l-5} \oplus (\lambda_{l-1}+\lambda_l) & \text{for
$i=l-3$}\\
\lambda_{l-4} \oplus 2\lambda_l \oplus 2\lambda_{l-1} & \text{for
$i=l-2.$}
\end{cases}$$

Since $2\lambda_1$ is in the range of $k  \ell$ and hence
invertible, it follows $2\lambda_1 \trunc 
\theta=\theta.$ Since $\lambda_1 \trunc \lambda_1$ contains
$2\lambda_1$ by Lemma \ref{lm:exist}, we conclude $2\lambda_1 \trunc
\lambda_1=\lambda_1$ 
so inductively $2\lambda_1 \trunc \lambda_{i}=\lambda_{i}$ for $i<l-1.$   Similarly
\begin{align*}
2\lambda_l \trunc \theta &= \lambda_{l-2},\\
2\lambda_{l-1} \trunc \theta &= \lambda_{l-2},\\
2\lambda_l \trunc 2\lambda_1 &= 2\lambda_{l-1},\\
2\lambda_{l-1} \trunc 2\lambda_1 &= 2\lambda_{l}.
\end{align*}

Since $\lambda_1 \trunc \theta$ contains $\lambda_1,$ it is clear that
$\lambda_1 \trunc \lambda_1$ contains $0,$ $\theta,$ and $2\lambda_1,$
each with multiplicity one.   Notice that by Equation (\ref{eq:AP}),
$\eta$ is not a direct summand of $\lambda \trunc \gamma$ if the
distance between $\lambda$ and $\eta$ is more than $||\gamma||$ (this
is argued explicitly in the proof of Lemma 2 of
\cite{Sawin02a}).  Noting that $||\lambda_1||=1$ and computing
$||\lambda-\lambda_1||$ for $\lambda=\lambda_i$ with $2<i<l-1$ and
$\lambda=\lambda_{l-1}, \lambda_l, (\lambda_{l-1}+\lambda_l)$ we
conclude
$$\lambda_1 \trunc \lambda_1= 0 \oplus \theta \oplus 2\lambda_1.$$
It then follows recursively that
\begin{align*}
\lambda_i \trunc \lambda_1 &=
\begin{cases}
\lambda_{i-1} \oplus \lambda_{i+1} & \text{for $1<i<l-2$}\\
\lambda_{l-3} \oplus (\lambda_{l-1} + \lambda_l) & \text{for
$i=l-2,$}
\end{cases}\\
(\lambda_{l-1} +\lambda_l) \trunc \lambda_1&=\lambda_{l-2} \oplus
2\lambda_{l-1} \oplus 2\lambda_l,\\
2\lambda_{l-1} \trunc \lambda_1 &=\lambda_{l-1} + \lambda_l,\\
2\lambda_{l} \trunc \lambda_1 &=\lambda_{l-1} + \lambda_l.
\end{align*}
Finally, we get recursively from this
\begin{align*}
\lambda_i \trunc \lambda_j &=
\begin{cases}
\lambda_{i-j} \oplus \lambda_{i+j} & \text{for $j<i<l-1$ and
$i+j<l-1$}\\
\lambda_{i-j} \oplus (\lambda_{l-1} +\lambda_l) & \text{for $j<i<l-1$ 
and 
$i+j= l-1, l+1$}\\
\lambda_{i-j} \oplus 2\lambda_{l-1} \oplus 2\lambda_l  & \text{for 
$j<i<l-1$ and
$i+j= l$}\\
\lambda_{i-j} \oplus \lambda_{2l-i-j} & \text{for $l>i>j$ and $i+j>l+1$}\\
0 \oplus 2\lambda_1 \oplus \lambda_{i+j} & \text{for $j=i<l-1$ and
$i+j<l-1$}\\
0 \oplus 2\lambda_1 \oplus (\lambda_{l-1} + \lambda_l)  & \text{for 
$j=i<l-1$ and
$2i= l-1, l+1$}\\
0 \oplus 2\lambda_1 \oplus 2\lambda_{l-1} \oplus 2\lambda_l  & 
\text{for $j=i<l-1$ and
$2i= l$}\\
0 \oplus 2\lambda_1 \oplus \lambda_{2l-2i} & \text{for $j=i<l$ and $2i>l+1,$}
\end{cases}\\
\lambda \trunc \lambda_j &=
\begin{cases}
\lambda_{l-1}+\lambda_l & \text{for $\lambda=2\lambda_{l-1}$ or $\lambda=
2\lambda_l$ and $j=1$}\\
\lambda_{l-j} & \text{for $\lambda=2\lambda_{l-1}$ or $\lambda=2\lambda_l$ and 
$1<j<l-1.$}
\end{cases}
\end{align*}

Thus if the closed subset $\Gamma$ contains $\lambda_1,$  $\lambda_2,$
$\lambda_{l-1},$ $\lambda_l$ or $\lambda_{l-1}+\lambda_l$ it
contains $\theta$ and is of the form $\Gamma_Z.$  If it contains only a
subset of $0,$ $2\lambda_1,$ $2\lambda_{l-1},$ and $2\lambda_l$ it is
of the form $\Delta_Z.$   If $\Gamma$ is not of the form $\Gamma_Z$ or
$\Delta_Z$ then it must contain $\lambda_j$ for some $2<j<l-1,$ so
suppose $j$ is the least such.  Then $\Gamma$ contains $\lambda_j,
\lambda_{2j}, \ldots, \lambda_{mj},$ where $m$ is the greatest such
that $mj<l-1.$  Again $\Gamma$ must contain every summand of
$\lambda_{mj} \trunc \lambda_j.$   Then  by the maximality of  $m$ we
know $(m+1)j>l-2,$ and since $\Gamma$ cannot contain $\lambda_{l-1} +
\lambda_l$ we know $(m+1)j \neq l-1,l+1.$  If $(m+1)j=l,$ then
$\Gamma$ contains $\{0,2\lambda_1, \lambda_j, \lambda_{2j}, \ldots,
\lambda_{mj},  2\lambda_{l-1}, 2\lambda_l\}.$  If it contained any
other $\lambda_i$ for $2<i<l-1$ then there would be $p$ with 
$|i-pj|<j,$ so $\lambda_{|i-pj|}$ would be in the set, contradicting 
the
minimality of $j.$  Since $\Gamma$ cannot contain any other weights in the
alcove, we conclude $\Gamma$ is of the form of set (5)
 in the theorem, with $m=n-1.$

On the other hand if $(m+1)j \neq l, $ the $(m+1)j>l+1,$ and therefore
$\lambda_{2l+1-(m+1)j}.$  Of course $2l+1-(m+1)j$ is a distance less
than $j$ from $mj,$ so if the difference is nonzero then again we contradict 
the minimality  of $j.$   Therefore the
distance between them is zero, so $2l+1=(2m+1)j,$ and $\Gamma$
contains $\{0, 2\lambda_1, \lambda_j, \cdots, \lambda_{mj}\}.$   Again
it cannot contain any other weight in the Weyl alcove without
contradicting the minimality of $j$  (if it contained $2\lambda_l$ or
$2\lambda_{l-1}$ it would contain $\lambda_{l-j},$  which is distinct
from $\lambda_{mj}$  but $l-j$ is less than $j$ away from $mj$) so
$\Gamma$ is set (6) in the theorem, with $n=2m+1.$

\item For $E_6,$ the weights of the Weyl alcove are $0,$ $\lambda_1,$
$2\lambda_1,$ $\lambda_2,$ $\lambda_3,$ $\lambda_5,$ $\lambda_6,$
$2\lambda_6,$ and $\lambda_1 + \lambda_6.$  A  closed subset
containing a  dull corner must
contain a nontrivial weight of the form $\theta+ \alpha,$ but the only
such weight is $\lambda_1+ \lambda_6,$ which is not a corner.

\item For $E_7,$ the weights in the alcove are  those in $k  \ell$  ($0$ and
$2\lambda_7$), the other corners ($\lambda_1=\theta,$ $\lambda_2$
and $\lambda_6$)  and one other ($\lambda_7$).  We have
\begin{align*}
2\lambda_7 \trunc \theta &= \lambda_6,\\
\theta  \trunc \theta &= 0 \oplus \lambda_6,\\
\lambda_6  \trunc \theta &= \theta \oplus 2\lambda_7,\\
\lambda_2  \trunc \theta &= \lambda_7,\\
\lambda_7  \trunc \theta &= \lambda_2 \oplus \lambda_7,
\end{align*} 
so 
\begin{align*}
\lambda_6 \trunc \lambda_6 &= 0 \oplus \lambda_6,\\
\lambda_2 \trunc \lambda_6&= \lambda_7,\\
\lambda_7 \trunc \lambda_6 &= \lambda_2 \oplus \lambda_7.
\end{align*}
Since every weight in $E_7$ is self-dual, $\lambda_7 \trunc \lambda_7$
consists of weights in the root lattice, and since
$N_{\lambda,\gamma}^\delta= N_{\lambda,\delta^*}^{\gamma^*}$ we can
read off from the previous equations
\begin{align*}
\lambda_7 \trunc \lambda_7 &= 0 \oplus \theta \oplus 2\lambda_7 \oplus
\lambda_6,\,\, \text{ so}\\
\lambda_2 \trunc \lambda_7 &= \theta \oplus \lambda_6\,\, \text{ and}\\
\lambda_2 \trunc \lambda_2 &= 0 \oplus 2\lambda_7.
\end{align*}

Thus the smallest closed subset containing $\lambda_6$ is
$\{0,\lambda_6\},$ containing $\lambda_2$ is
$\{0,\lambda_2,2\lambda_7\},$ and containing any other nonunit is a
subset containing $\theta.$

\item For $E_8,$ the only weights are $0,$ $\lambda_1$ and
$\lambda_8=\theta.$  As above we see that $\theta \trunc \theta= 0 \oplus
\lambda_1,$ $\lambda_1 \trunc \theta= \theta,$ and hence $\lambda_1
\trunc \lambda_1=0.$  Thus $\lambda_1$ is invertible,
$\{\lambda_1,0\}$ is a closed subset, and every other closed
subset contains $\theta.$

\item For $F_4$ and $G_2$ there is nothing to prove since the only
elements are $0,$ $\theta,$ $\beta,$ and $2\beta,$ where $2\beta$ is
a short weight.
\end{itemize}

\subsection{If $\mathbf{k=3}$}    Since now for all dull corners we know $\lambda \trunc \lambda^\dagger$ contains something
nontrivial in $\beta \trunc \beta,$  we may assume our closed set
contains a  corner which is of the form $\beta+\alpha$  for $\alpha$ short.

 In the nonsimply-laced case there are no such
corners, so we need only consider the simply-laced 
case. The only corners are
fundamental weights with $(\lambda_i,\theta)=3$ together with the range of $k
 \ell.$

\begin{itemize}
\item For $D_l$ there are no fundamental weights  with
$(\lambda_i,\theta)=3,$ and the range of $k\ell,$ ($3\lambda_1,$
$3\lambda_{l-1},$ and $3\lambda_l$) contains nothing in the root lattice.

\item For $E_6,$ neither $3\lambda_1$ nor $3\lambda_6$ is the sum of
two roots, so only $\lambda_4$ is such a  corner.  Since
$\lambda_4 \trunc \theta$ contains at least three summands, $\lambda_4
\trunc \lambda_4$ contains two distinct nontrivial summands of $\theta
\trunc \theta,$ so it must contain a noncorner.

\item For $E_7,$ since $3\lambda_7$ is not in the root lattice, only
$\lambda_3$ and $\lambda_5$ are such  corners, 
and  only $\lambda_3$ among all corners is of the form $\theta +
\alpha.$  Thus  any closed subset containing a corner contains
$\lambda_3,$ so since $\lambda_3 \trunc \theta$ contains at least
three summands, $\lambda_3 \trunc \lambda_3$ contains two distinct
nontrivial summands of $\theta \trunc \theta,$ one of which must not
be a corner.

\item For $E_8,$ Only $\lambda_2$ and $\lambda_7$ are corners, and
neither is of the form $\theta+ \alpha.$
\end{itemize}

\subsection{If $\mathbf{k=4}$}  Again we need consider only the
simply-laced case, and we need only consider corners of the form
$\theta +\alpha,$ which means the corner $2\theta.$  Now in each case
$\theta=\lambda_i$ for some $i,$ so in addition to the weights in the
chart, we have $2\theta -\alpha_i$ is a summand of $2\theta \trunc
\theta,$ so $2\theta \trunc 2\theta$ contains two distinct nontrivial
summands of $\theta \trunc \theta,$ one of which must not be a corner.

\subsection{If $\mathbf{k>4}$}  For any dull corner $\lambda$  $\lambda \trunc \lambda^*$ contains a nontrivial summand of
$\theta \trunc \theta,$ which cannot be a corner.
\end{proof}

\section{Modular Categories and Closed Subsets of the Weyl Alcove}

As alluded to in the introduction, there is a well-defined procedure
for  constructing  a modular or spin modular category out a ribbon
$*$-category (all of our ribbon categories inherit the $*$-structure
from the original ribbon category).  These techniques were developed by
M\"uger \cite{Muger99} and Brugui\`eres \cite{Bruguieres00}.

First one should identify the subcategory of degenerate objects, which
are simple objects (in our case irreducible representations) $\lambda$
such that $R_{\lambda,\gamma}=R^{-1}_{\gamma,\lambda}$ for every
simple object $\gamma$ in the ribbon category, where $R$ is the
$R$-morphism associated to a crossing.  These come in two sorts, even or
odd, according to whether the effect of the full twist $C_\lambda$ is
multiplication by one or minus one.  If all are even, one can quotient
by them to get a modular category.  If in addition they are all
invertible and form a cyclic group, \cite{Sawin02b} offers a detailed description of the
invariant and TQFT in  terms of the original category. If there are
odd degenerate objects, one can still quotient by the even degenerate
objects and the result is a spin-modular category which gives an
invariant of spin three-manifolds \cite{Sawin02b}.

For the ribbon categories associated to the closed sets of the form
$\Gamma_Z$ and $\Delta_Z$ \cite{Sawin02a}  proves  that all
degenerate objects are invertible and identifies in which cases they
are all odd. Below we identify for each exceptional closed  category 
what the degenerate objects are, when they are even, when they are 
invertible, and when the resulting TQFT is equivalent to a 
nonexceptional TQFT.

\subsection{The case $\mathbf{E_7}$ with $\mathbf{\{0,\lambda_6\}}$}
Here $\lambda_6 \trunc \lambda_6=0 \oplus \lambda_6.$

For $E_7$ we have $h=18$ so $k+h=20.$  We have
$(\lambda_6,\lambda_6)=4$ and $(\lambda_6,\rho)=26,$ so
$$C_{\lambda_6}=e^{2\pi i (4+ 2 \cdot 26)/(2\cdot 20)}=e^{4\pi
i/5}.$$
Since this is not $\pm 1$ then $\lambda_6$ is not degenerate, so the
set is in itself modular.  

Straightforward calculations show that this ribbon category is
determined up to isomorphism by the link invariant, which is in turn
determined by a skein relation, the same skein relation that
determines the $SO(3)$ theory at $k=3,$ and thus this ribbon category
is isomorphic to the $SO(3)$ level $3$ category.

\subsection{The case $\mathbf{E_7}$ with $\mathbf{\{0,\lambda_2,2\lambda_7\}}$}
Here $\lambda_2 \trunc \lambda_2=0 \oplus 2\lambda_7$ and $\lambda_2
\trunc 2\lambda_7=\lambda_2.$

Again $h+k=20,$ $(\lambda_2,\lambda_2)=7/2,$  $(\lambda_2,\rho)=49/2,$
$(2\lambda_7,2\lambda_7)= 6$ and $(2\lambda_7,\rho)=27$ so
\begin{align*}
C_{\lambda_2}&=e^{2\pi i (7/2 + 2 \cdot 49/2)/(2 \cdot 20)}=e^{5\pi i
/8}\\
C_{2\lambda_7}&=e^{2 \pi i (6+ 2 \cdot 27)/(2\cdot 20)}=-1.
\end{align*}
Of course
$$R_{\lambda_2,2\lambda_7}R_{2\lambda_7,\lambda_2}=C_{\lambda_2}
C^{-1}_{\lambda_2} C^{-1}_{2\lambda_7}=-1$$
so neither $2\lambda_7$ nor $\lambda_2$ is degenerate and the theory is
modular.  Again direct calculation shows that this theory is 
isomorphic to the $\mathrm{SU}(2)$ theory at level $k=2$ but this time with a 
nonstandard choice for $q^{1/4}=e^{13\pi i/8}$ (see \cite{Sawin05b}).

\subsection{The case $\mathbf{E_8}$ with $\mathbf{\{0,\lambda_1\}}$}
Here $\lambda_1 \trunc \lambda_1=0.$

For $E_8$ $h=30$ so $k+h=32.$  $(\lambda_1,\lambda_1)= 4$ and
$(\lambda_1,\rho)=46$ so
$$C_{\lambda_1}=e^{2\pi i (4 + 2\cdot 46)/(2\cdot 32)}=-1.$$
This is an odd degenerate object and of course the ribbon category is
isomorphic to that for $SO(3)$ at $k=2.$  As argued in
\cite{Sawin02b}, this gives the trivial invariant of spin
three-manifolds.

\subsection{The case $\mathbf{B_l}$ with $\mathbf{\{0,2\lambda_1, \lambda_j,
\lambda_{2j}, \ldots, \lambda_{(n-1)j/2}\}}$ where $\mathbf{2l+1=nj}$}
Here $h=2l-1$ so $k+h=2l+1.$  Also, $(\rho, \lambda_i)=li-i^2/2$ and
$(\lambda_i,\lambda_i) =i.$  Thus
$$C_{2\lambda_1}= e^{\pi i (4 + 2(2l-1))/(2l+1)}=1$$
and thus
$$R_{2\lambda_1,\lambda_{mj}}R^{-1}_{\lambda_{mj},2\lambda_1}=1$$
so $2\lambda_1$ is even degenerate.   On the other hand
$$C_{\lambda_{mj}}=e^{\pi i (mj-m^2j^2/(2l+1))}$$
so 
\begin{align*}
C_{\lambda_{(m \pm p)j}} C_{\lambda_{mj}}^{-1} C_{\lambda_{pj}}^{-1}&=
e^{\mp 2 \pi i mpj^2/(2l+1)},\\
C_{\lambda_{2l+1-(m + p)j}} C_{\lambda_{mj}}^{-1} C_{\lambda_{pj}}^{-1}&=
e^{- 2 \pi i mpj^2/(2l+1)},\\
C_{2\lambda_1}C_{\lambda_{mj}}^{-1} C_{\lambda_{mj}}^{-1}= 
C_{0} C_{\lambda_{mj}}^{-1}  C_{\lambda_{mj}}^{-1}&=e^{2\pi i 
m^2j^2/(2l+1)}.
\end{align*}
This will be  $1$ for all $p,$ and
thus $\lambda_{mj}$ will be degenerate,  exactly if $m$ is a multiple
of $n/d,$ where $d$ is the greatest  common divisor of $n$ and $j.$ 
If $m=rn/d$ then 
$$C_{\lambda_{mj}}=e^{\pi i (mj-m^2j^2/(2l+1))}=e^{\pi i(r
(2l+1)/d-r^2(2l+1)/d^2)}.$$
Since $(2l+1)/d$ and $(2l+1)/d^2$ are both odd, this is one whether
$r$ is even or odd, and thus all such $\lambda_{mj}$ are even.  Thus we
get a modular category which is a quotient by the $\Z/2$ action if $m$
and $j$ are relatively prime, and by the set of even simple degenerates 
$$\{0,
2\lambda_1, \lambda_{(2l+1)/d}, \lambda_{2(2l+1)/d}, \ldots,
\lambda_{(d-1)(2l+1)/(2d)}\}$$
if $(n,j)=d \neq 1.$  Notice that when $d \neq 1$  the set of simple
degenerates does not form a group, which is to 
say that there are noninvertible degenerate objects. In fact one can 
check that the subcategory generated by these representations is 
isomorphic to the 
representation theory of the nonabelian group presented by
$\bracket{x,y\, |\, 
x^2=y^2=(xy)^d=1}$  (any eigenvector of $xy$ generates a
two-dimensional irreducible subrepresentation if the eigenvalue of
$xy$ is 
a nontrivial $d$th root of unity and a one-dimensional
subrepresentation if the eigenvalue is $1,$ recalling that $d$ is
necessarily odd.  These two-dimensional representations are classified
by the eigenvalue of $xy,$ and the one-dimensional by the eigenvalue
of $x.$) For example, when $l=13$ and $j=3$, we have $n=9,$
$d=3,$ and the subcategory generated by 
\[\{0, 2\lambda_1, \lambda_3, \lambda_6, \lambda_9, \lambda_{12}
\}\]
contains as its trivial subcategory 
\[\{0, 2\lambda_1,  \lambda_9
\}\]
where $2\lambda_2 \trunc 2\lambda_1=0,$ $2\lambda_1 \trunc
\lambda_9=\lambda_9,$ and $\lambda_9 \trunc \lambda_9= 0 \oplus
2\lambda_1 \oplus \lambda_9.$  

This is the
first example of which the author is aware of a nonsymmetric ribbon category with
noninvertible degenerate objects.  Whether the resulting quotient 
gives a truly new TQFT is not clear, and in any case this example is 
worthy of further study.

\subsection{The case $\mathbf{D_l}$ with $\mathbf{\{0,2\lambda_1, \lambda_j,
\lambda_{2j}, \ldots , \lambda_{(n-1)j}, 2\lambda_{l-1}, 2\lambda_l\}}$
where $\mathbf{l=nj}$}
Here $h=2l-2$ so $k+h=2l,$  $(\rho,\lambda_i)=i(2l-i-1)/2,$ and
$(\lambda_i,\lambda_i)=i,$ so 
\begin{align*}C_{2\lambda_1}&=e^{\pi i (4+2\cdot 2\cdot(2l-2)/2)/(2l)}=1,\\
C_{2\lambda_l}&= e^{\pi i(4l+ 2\cdot 2l(2l-l-1)/2)/(2l)}=e^{\pi i
(l+1)},\\
C_{2\lambda_{l-1}}&= e^{\pi i(4(l-1)+ 2\cdot 2(l-1)(2l-l)/2)/(2l)}=e^{\pi i
(l+1)},\\
C_{\lambda_{mj}}&=e^{\pi i (mj+2mj(2l-mj-1)/2)/(2l)}=e^{\pi 
i(mj-m^2j^2/(2l))},
\end{align*}
so $2\lambda_1$ is even degenerate and
\begin{align*}
R_{2\lambda_l,\lambda_{mj}}R^{-1}_{\lambda_{mj},2\lambda_l}&=C_{\lambda_{l-mj}}
C_{\lambda_{mj}}^{-1} C_{2\lambda_l}^{-1}=e^{\pi i(-l/2-mj-1)},\\
R_{2\lambda_{l-1}l,\lambda_{mj}}R^{-1}_{\lambda_{mj},2\lambda_{l-1}}&=
C_{\lambda_{l-mj}} C_{\lambda_{mj}}^{-1} C_{2\lambda_{l-1}}^{-1}=e^{\pi
i(-l/2-mj-1)},
\end{align*}
and $2\lambda_l,$ $2\lambda_{l-1}$  are thus  degenerate exactly when 
$j$ is even and $l/2$ is odd.  They are always odd degenerate.    Also 
\begin{align*}
    C_{\lambda_{(m \pm p)j}}C_{\lambda_{pj}}^{-1} 
    C_{\lambda_{mj}}^{-1} &=
     e^{\mp \pi i mpj^2/l},\\
     C_{\lambda_{l-(m + p)j}} C_{\lambda_{mj}}^{-1} C_{\lambda_{pj}}^{-1}&=
     e^{-  \pi i mpj^2/l},\\
     C_{2\lambda_1}C_{\lambda_{mj}}^{-1} C_{\lambda_{mj}}^{-1} = 
    C_{0}  C_{\lambda_{mj}}^{-1}   C_{\lambda_{mj}}^{-1} &= e^{-\pi i 
    m^2j^2/l},\\
    C_{2\lambda_{l-1}} C_{\lambda_{mj}}^{-1} C_{\lambda_{l-mj}}^{-1} =  
    C_{2\lambda_{l}} C_{\lambda_{mj}}^{-1} 
    C_{\lambda_{l-mj}}^{-1}     &=  e^{\pi i (1+mj^2/l -mj + l/2)}.
\end{align*}
  This will be $1$ for all $p,$ and thus $\lambda_{mj}$ will be 
  degenerate, exactly when $l$ is even, $l/2$ is odd, and either $m$
  is an even  multiple 
  of   $n/d,$ where $d$ 
  is the greatest common divisor of $n$ and $j,$  or $j$ is even and
  $m$ is a multiple of $l/d.$ 
   If $m=rn/d,$ then 
  $$C_{\lambda_{mj}}=e^{\pi i (rl/d-r^2l/(2d^2))}$$
  which is  $1$ if and only if   $r$ is even.  

Thus the set of 
  simple degenerate elements consists of
  $$\{0,2\lambda_1\}$$
  if  (a) $l$ is odd, or (b) $l$ and $l/2$ are even, or (c) $l$ is even, $l/2$ is 
  odd and $\gcd(n,j)=1.$ In this case the quotient is modular. 

The set 
  of simple degenerate objects consists of  
  $$\{0, 2\lambda_1, \lambda_{2l/d},\lambda_{4l/d}, \ldots, 
  \lambda_{(d-1)l/d}\}$$
  if $l$ is even, $j$ and $l/2$ are odd and $\gcd(n,j) \neq 1.$  In 
  this case all the degenerate objects are even and the quotient is
  modular.  The subcategory of  degenerate  
  objects is isomorphic to the category of representations of the 
  group  $\bracket{x,y\, |\, 
x^2=y^2=(xy)^{(d+1)/2}=1}.$  

The set of simple degenerate objects is 
  $$\{0, 2\lambda_1, \lambda_{l/d},\lambda_{2l/d}, \ldots, 
  \lambda_{(d-1)l/d},2\lambda_{l-1},2\lambda_l\}$$
  if $l$ and $j$ are even and $l/2$ is odd (here the set is understood to have 
  just four elements if $d=1$).  In this  case the odd multiples of
  $l/d$ and the  
  last two entries give  odd degenerate objects.   Notice  the 
  even degenerate objects are isomorphic to the category of 
  representations of the  group  $\bracket{x,y\,  |\, 
x^2=y^2=(xy)^{d}=1}.$ 

\subsection{The case $\mathbf{D_l}$ with $\mathbf{\{0,2\lambda_1, \lambda_j,
\lambda_{2j}, \ldots , \lambda_{(n-1)j/2}\}}$
where $\mathbf{2l=nj}$ for $\mathbf{n}$ odd}  
Again 
\begin{align*}C_{2\lambda_1}&=e^{\pi i (4+2\cdot 2\cdot(2l-2)/2)/(2l)}=1,\\
C_{\lambda_{mj}}&=e^{\pi i (mj+2mj(2l-mj-1)/2)/(2l)}=e^{\pi 
i(mj-m^2j^2/(2l))},
\end{align*}
so $2\lambda_1$ is still even degenerate and 
\begin{align*}
    C_{\lambda_{(m \pm p)j}}C_{\lambda_{pj}}^{-1} C_{\lambda_{mj}}^{-1}&=
e^{\mp \pi i mpj^2/l}\\
     C_{\lambda_{l-(m + p)j}} C_{\lambda_{mj}}^{-1} C_{\lambda_{pj}}^{-1}&=
     e^{-  \pi i mpj^2/l},\\
     C_{2\lambda_1}C_{\lambda_{mj}}^{-1} C_{\lambda_{mj}}^{-1} = 
    C_{0}  C_{\lambda_{mj}}^{-1}   C_{\lambda_{mj}}^{-1} &= e^{-\pi i 
    m^2j^2/l}.
    \end{align*}
  This will be $1$ for all $p,$ and thus $\lambda_{mj}$ will be 
  degenerate, exactly when $m$ is a multiple of $n/d,$ where $d$ 
  is the greatest common divisor of $n$ and $j.$   If $m=rn/d,$ then 
  $$C_{\lambda_{mj}}=e^{\pi i (2rl/d-2r^2l/d^2)}$$
  which is always $1.$   Thus the set of simple degenerate objects is
  $$\{0,2\lambda_1, \lambda_{l/d},\lambda_{2l/d}, \ldots, 
  \lambda_{(d-1)l/d}\},$$
all degenerate objects are even, and again the category of even 
degenerate objects is isomorphic to the representation category of 
the group $\bracket{x,y \, |\, 
x^2=y^2=(xy)^d=1}.$

\section{Tensor Closed Implies Closed} \label{se:tensor}

The definition of closed involves two conditions:  closed under the truncated tensor
product  and closed under duality.  Under the sort of
conditions found in the quantum group examples the first condition
actually implies the second.    In
particular for any sum of weights in the Weyl alcove,  the set of
weights appearing as summands of  tensor powers
of that sum is one of the closed subsets classified above.   This is of
particular relevance to skein-theoretic and Young diagrammatic approaches to the link
invariants (see, e.g., Turaev and Wenzl \cite{TW93})  where all the link information is
recovered from cabling  (i.e.  tensor powers) of an invariant corresponding to
one particular weight   (corresponding to the fundamental representation).  The
approach to  the proof of the following proposition
was suggested to the author by A. Liakhovskaia.  

\begin{prop} \label{pr:tensor}
If $\lambda$ is in the Weyl alcove then there
is an $n$ such that $\lambda^{\trunc n}$ contains $\lambda^\dagger$ as a summand.
\end{prop}

\begin{proof} We shall actually show there is an $m$ such that $\lambda^{\trunc m}$
contains the weight $0$ as a summand:  Of course $n=m-1$ then suffices for the
proposition.  

For any two elements of the Weyl alcove $\lambda$ and $\gamma$ let 
$S_{\lambda,\gamma}$ be
the value of the link invariant on the Hopf link with its two components
 labeled by
$\lambda$ and $\gamma$ respectively.  Recall from \cite{Sawin02a} that viewed
 as a
$|\Lambda_0|$-by-$|\Lambda_0|$ matrix $S$ is nondegenerate,  and that 
\[
\sum_{\gamma \in \Lambda_0} \qdim(\gamma) S_{\lambda,\gamma}=
\delta_{\lambda,0}
\sum_{\gamma \in \Lambda_0} \qdim(\gamma)^2.
\]
Thus if $\sum_\gamma S_{\lambda^{\tensor m},\gamma}$ is nonzero
$\lambda^{\tensor m}$ contains $0$ as a summand.  Thus it suffices to
show $\sum_\gamma S_{\lambda^{\tensor m},\gamma}$ is nonzero for some
$m.$ Dividing by $\qdim(\lambda)^m$ we see
\begin{equation} \label{eq:rep_sum}
\sum_\gamma S_{\lambda^{\trunc m},\gamma}\qdim(\gamma)/\qdim(\lambda)^m=\sum_\gamma
\left[S_{\lambda,\gamma}/\qdim(\gamma)\qdim(\lambda)\right]^m\qdim(\gamma)^2.
\end{equation}

Now 
\begin{gather*}
S_{\lambda,\gamma}=\sum_{\mu \in \Lambda_0} N_{\lambda,\gamma}^\mu C_\mu
C_\lambda^{-1} C_\gamma^{-1} \qdim(\mu),\\
S_{\lambda \trunc \gamma,\mu}=S_{\lambda,\mu} 
S_{\gamma,\mu}/\qdim(\mu),\,\,\text{ and}\\
\end{gather*}
 so
since $\qdim(\gamma)$ is positive, $C_\mu$ is a root of unity, and 
$$\qdim(\lambda) \qdim(\gamma)= \sum_{\mu \in \Lambda_0} N_{\lambda,\gamma}^\mu  
\qdim(\mu),$$
we see that $|S_{\lambda,\gamma}/(\qdim(\lambda) \qdim(\gamma))| \leq
1,$ and $S_{\lambda,\gamma}/(\qdim(\lambda)\qdim(\gamma))$ is a root
of unity when the absolute value equals one.  
 
  The quantity in square brackets on the right-hand side of Equation
  (\ref{eq:rep_sum}) has modulus $\leq 1$ for
all values of $\gamma,$ and for at least one term in the sum ($\gamma=0$) has modulus
equal to $1.$  So for very large $m$ this sum is dominated by  terms where
the modulus is equal to $1.$  In each of these terms the quantity
$S_{\lambda,\gamma}/\qdim(\gamma)\qdim(\lambda)$ is a root of unity, so for infinitely
many values of $m$ the value of $(S_{\lambda,\gamma}/\qdim(\gamma)\qdim(\lambda))^m$ is
equal to $1$ simultaneously for all of the values of $\gamma$ for which the ratio has
modulus
$1.$  Thus for sufficiently large $m$ the sum must be positive.
\end{proof}


\def\cprime{$'$}

\end{document}